\documentclass[12pt, a4paper]{article}

\usepackage{amsmath,amsfonts}
\usepackage{mathrsfs} %to get script letters \mathscr{L} yields L
\usepackage{parskip}
\usepackage[T1]{fontenc}
\usepackage[utf8]{inputenc}
\usepackage{authblk}
\usepackage{enumerate}

\newtheorem{theorem}{Theorem}
\newtheorem{lemma}{Lemma}

\def \Liminf{\mathop{\underline{\lim}}\limits}

\def\Pb{\mathbf{P}}

\def\Ex{\mathbf{E}}

\def\KK{\mathbb{K}}

\def\II{\mathbb{I}}

\def\UU{\mathbb{U}}
\def\VV{\mathbb{V}}

\def\sgn{{\rm sgn}} 
\def\1{\mbox{1\hspace{-.25em}I}}
\begin{document}
\title{On  Parameter Estimation for Cusp-type Signals}
\author[1]{O.V. Chernoyarov}
\author[2]{S. Yu. Dachian}
\author[3]{Yu.A. Kutoyants}
\affil[1,3]{National Research University ``MPEI'', Moscow, Russia, }
\affil[2]{University of Lille 1, Lille,  France,}
\affil[3]{University of  Maine,  Le Mans,  France}

\date{}
\maketitle

\begin{abstract}
We consider the problem of parameter estimation by the observations of
deterministic signal in white gaussian noise. It is supposed that the signal
has a singularity of {\it cusp}-type. The properties of the maximum likelihood
and bayesian estimators are described in the asymptotics of small
noise. Special attention is paid to the problem of parameter estimation in the
situation of misspecification in regularity, i.e.; the statistician supposes
that the observed signal has this singularity, but the real signal is
smooth. The rate and the asymptotic distribution of the maximum likelihood
estimator in this situation are described. 
\end{abstract}
\noindent MSC 2000 Classification: 62M02,  62G10, 62G20.

\noindent {\sl Key words}: \textsl{Parameter estimation, cusp-type
  singularity, small noise, misspecification.}  

\section{Introduction}

Consider the problem of parameter estimation by the observations
$X^T=\left(X_t,0\leq t\leq T\right)$ of the 
signals in White Gaussian Noise (WGN)
\begin{align}
\label{01}
{\rm d}X_t=S\left(\vartheta ,t\right){\rm d}t+\varepsilon {\rm d}W_t,\quad
X_0=0,\quad 0\leq t\leq T.
\end{align}
Here $S\left(\vartheta ,t\right)$ is a known function (signal), $W_t,0\leq
t\leq T$ is a Wiener process and $\vartheta \in \Theta =\left(\alpha ,\beta
\right)$ is unknown parameter. 

We have to estimate the parameter $\vartheta $ by continuous time observations
$X^T$ and to describe the properties
of estimators in the asymptotics of {\it small noise}, i.e., the parameter
$\varepsilon \in (0,1]$ is known and the asymptotics corresponds to  $\varepsilon \rightarrow
0$. 

It is known that if the signal $S\left(\vartheta ,\cdot \right)$ is a smooth
function of $\vartheta $, then the maximum likelihood estimator and bayesian
estimators are consistent, asymptotically normal
\begin{align*}
\varepsilon ^{-1}\left(\hat\vartheta _\varepsilon-\vartheta
\right)\Longrightarrow {\cal N}\left(0,\II\left(\vartheta
\right)^{-1}\right),\qquad \varepsilon ^{-1}\left(\tilde\vartheta
_\varepsilon-\vartheta 
\right)\Longrightarrow {\cal N}\left(0,\II\left(\vartheta \right)^{-1}\right),
\end{align*}
we have the convergence of all polynomial moments
and the both estimators are asymptotically efficient \cite{IH1}. Here
$\II\left(\vartheta 
\right)$ is the Fisher information 
\begin{align}
\label{0}
\II\left(\vartheta \right)=\int_{0}^{T}\dot S\left(\vartheta ,t\right)^2\,{\rm
  d}t. 
\end{align}

Here and in the sequel dor means derivation w.r.t. $\vartheta $. If the signal
$S\left(\vartheta ,t\right)=S\left(t-\vartheta \right)$, where
$S\left(t\right)$ is a discontinuous function of $t$, say, has a jump at the
point $t=0$. Then $\II\left(\vartheta \right)=\infty $, the MLE $\hat\vartheta
_\varepsilon$ and BE $\tilde\vartheta _\varepsilon$ have the rate of
convergece $\varepsilon ^2$ with different limit distributions:
\begin{align*}
\varepsilon ^{-2}\left(\hat\vartheta _\varepsilon-\vartheta
\right)\Longrightarrow \hat u,\qquad \quad \varepsilon
^{-2}\left(\tilde\vartheta _\varepsilon-\vartheta 
\right)\Longrightarrow \tilde u ,
\end{align*}
and asymptotically efficient are bayesian estimators only. Here $\Ex
\left(\hat u\right)^2> \Ex \left(\tilde u\right)^2$. For the proofs see
\cite{IH2}.

We are interested by the properties of the MLE $\hat\vartheta _\varepsilon $ in
the case of observations \eqref{01}, where the signal $S\left(\vartheta
,t\right)$ has a singularity of the {\it cusp}-type, i.e.; at the vicinity of
the point $t=\vartheta $ it  has the representation $S\left(\vartheta
,t\right)\approx a\left|t-\vartheta \right|^\kappa  $, where $\kappa \in
(0,\frac{1}{2})$. Note that for these values of $\kappa $   we have
$\II\left(\vartheta \right)=\infty $. 

The problem of parameter estimation for cusp-type singular density function
by i.i.d. observations 
was considered in \cite{PR68}. It was shown that the MLE $\hat\vartheta _n$
has limit distribution with the rate 
\begin{align*}
n^{\frac{1}{2\kappa +1}} \left(\hat\vartheta _n-\vartheta
\right)\Longrightarrow \hat \eta  .
\end{align*}
The exhaustive study of singular estimation problems for i.i.d. observations
including cusp-type singularity can be found in \cite{IH81}. For stochastic
processes observed in continuous time the similar problems were considered in
\cite{D03} for inhomogeneous Poisson processes and in \cite{DK03} for ergodic
diffusion processes. 

This work is devoted to two  problems. The first one is to describe the
asymptotics of the MLE and BE in the case of signal with cusp-type
singularity. It is shown that 
\begin{align*}
\varepsilon ^{-\frac{2}{2\kappa +1} }\left(\hat\vartheta _\varepsilon -\vartheta
\right)\Longrightarrow \hat \xi ,\qquad \quad \varepsilon ^{-\frac{2}{2\kappa
    +1} } \left(\tilde\vartheta _\varepsilon -\vartheta \right)\Longrightarrow
\tilde \xi 
\end{align*}
where $\hat \xi $ and $\tilde \xi $ are two different r.v.'s, $ \Ex (\hat
\xi ^2)>\Ex (\tilde  \xi ^2)$.  The second
problem is to study the properties of the MLE, when the signal supposed by the
statistician (theoretical) has cusp-type singularity, but the real signal is
smooth (regular). We show that
\begin{align*}
\varepsilon ^{-\frac{2}{3-2\kappa } }\left(\hat\vartheta _\varepsilon -\hat\vartheta
\right)\Longrightarrow \hat \zeta .
\end{align*} 
Here $\hat\vartheta$ is the value of $\theta $ which minimizes the
corresponding Kulback-Leibler distance. The proofs are carried out following
two general results by Ibragimov and Khasminskii (Theorems 1.10.1 and 1.10.2
in \cite{IH81}), i.e., we verify the conditions of these theorems for our
model of observations. 

Note that the similar problem of misspecification was considered in the work
\cite{CKT15}, where the  signal chosen by the statistician (theoretical model)
has discontinuity, but the real signal is smooth. It is shown that 
\begin{align*}
\varepsilon ^{-\frac{2}{3 } }\left(\hat\vartheta _\varepsilon -\hat\vartheta
\right)\Longrightarrow \hat \eta .
\end{align*} 
We discuss as well the problem of estimation $\kappa $.   The presented work
is a contuniation of the study \cite{CKT15}.

\section{Main  result}

Let us consider the problem of parameter estimation by the observations (in
continuous time) of the deterministic signal in the presence of White Gaussian
Noise (WGN) of small intensity
\begin{align}
\label{0111}
{\rm d}X_t=S\left(\vartheta_0 ,t\right){\rm d}t+\varepsilon {\rm d}W_t,\quad
X_0=0,\quad 0\leq t\leq T, 
\end{align}
where the unknown parametr $\vartheta_0 \in \Theta =\left(\alpha ,\beta
\right)$.  We are interested
by the  behavior of the estimators of this parameter  in the asymptotics of
{\it small noise}, i.e., as $\varepsilon \rightarrow 0$.

Suppose that the signal $S\left(\vartheta ,t\right)$ has {\it cusp}-type
singularity $$ S\left(\vartheta ,t\right)=a\left|t-\theta \right|^\kappa
+h\left(t,\vartheta \right), $$ where $0<\alpha<\vartheta <\beta <T $ and
$\kappa \in (0,\frac{1}{2})$. The function $h\left(\vartheta ,t\right)$ is
continuously differentiable w.r.t. $\vartheta $ and has bounded derivative. 

The likelihood ratio function is 
\begin{align*}
V\left(\vartheta ,X^T\right)=\exp\left\{\frac{1}{\varepsilon ^2}\int_{0}^{T}
S\left(\vartheta ,t\right){\rm d}X_t- \frac{1}{2\varepsilon ^2}\int_{0}^{T}
S\left(\vartheta ,t\right)^2{\rm d}t\right\},\quad \vartheta \in \Theta 
\end{align*}
(see \cite{LS01}) and the MLE $\hat\vartheta _\varepsilon $ is defined by the equation 
\begin{align}
\label{mle}
V\left(\hat\vartheta _\varepsilon ,X^T\right)=\sup_{\vartheta \in\Theta }V\left(\vartheta ,X^T\right).
\end{align} 

Suppose that $\vartheta $ is a random variable with continuous, positive
density function 
$p\left(\vartheta \right), \alpha <\vartheta <\beta $. The bayesian estimator
(BE) $\tilde\vartheta _\varepsilon $ with quadratic loss function is 
\begin{align}
\label{be}
\tilde\vartheta _\varepsilon=\frac{\int_{\alpha }^{\beta }\theta p\left(\theta
  \right)V\left(\theta ,X^T\right){\rm d}\theta }{\int_{\alpha }^{\beta }
  p\left(\theta \right)V\left(\theta ,X^T\right){\rm d}\theta}. 
\end{align}

We are interested by the properties of the estimators $\hat\vartheta
_\varepsilon $ and $\tilde\vartheta _\varepsilon $ in the asymptotics $\varepsilon
\rightarrow 0$.

 Note that the Fisher
information is not finite and we have a singular problem of parameter
estimation. Introduce the Hurst parameter $H=\kappa +\frac{1}{2}$ and
double-side fractional Brownian motion (fBm) $W^H\left(u\right), u\in R$.
%\begin{align*}
%W^H\left(u\right)=\begin{cases}
%W^H_+\left(u\right),& \text{if $u>0$},\\
%W^H_-\left(-u\right),& \text{if $u\leq 0$},
%\end{cases}
%\end{align*}
%Here $W^H_+\left(\cdot \right) $ and $W^H_-\left(\cdot \right) $ are two
%independent fractional Brownien movements. 
Recall, that $\Ex W^H_+\left(u
\right)=0 $ and 
\begin{align}
\label{fBm}
\Ex W^H_+\left(u \right)W^H_+\left(v
\right)=\frac{1}{2} \left[
  \left|u\right|^{2H}+\left|v\right|^{2H}-\left|u-v\right|^{2H}\right],\quad
u,v\in R. 
\end{align}

Introduce two  random variables $\hat\xi$
and $\tilde\xi$  by the relations 
\begin{align*}
Z(\hat\xi) =\sup_{u\in R}Z\left(u\right),\qquad \tilde \xi
=\frac{\int_{}^{}uZ\left(u\right){\rm d}u}{\int_{}^{}Z\left(u\right){\rm d}u},
\end{align*}
where the process 
\begin{align}
\label{z}
Z\left(u\right)=\exp \left\{\Gamma W^H\left(u\right)-\frac{\Gamma
   ^2}{2}\left|u\right|^ {2H}\right\},\qquad u\in R.
\end{align}
Here
\begin{align*}
\Gamma^2={a^2}\int_{-\infty }^{\infty }\left[ \left|v-
    1\right|^\kappa  -\left|v\right|^\kappa \right]^2{\rm d}v.
\end{align*}
Introduce as well the process
\begin{align*}
Z^o\left(v\right)=\exp \left\{ w^H\left(v\right)-\frac{1}{2}\left|v\right|^
{2H}\right\},\qquad v\in R 
\end{align*}
and the corresponding random variables $\hat\xi _o$ and $\tilde \xi _o$ by
the relations
\begin{align*}
Z(\hat\xi_o) =\sup_{v\in R}Z^o\left(v\right),\qquad \tilde \xi_o
=\frac{\int_{}^{}vZ^o\left(v\right){\rm d}v}{\int_{}^{}Z^o\left(v\right){\rm d}v}.
\end{align*}
Note that 
\begin{equation}
\label{equ}
\hat\xi =\frac{\hat\xi_o}{\Gamma ^{\frac{1}{H}}},\qquad \quad \tilde\xi
=\frac{\tilde\xi_o}{\Gamma ^{\frac{1}{H}}} .
\end{equation}
The proof of \eqref{equ} follows immediately from the change of variables $u=\Gamma
^{\frac{1}{H}}v$ in  $Z\left(u\right)$.

Adymptotically efficient estimators we define with the help of the following
lower bound. For all $\vartheta _0\in \Theta $ and all estimators
$\bar\vartheta _\varepsilon $ we have the relation
\begin{equation}
\label{lb}
\lim_{\delta \rightarrow 0}\Liminf_{\varepsilon \rightarrow 0}
\sup_{\left|\vartheta -\vartheta _0\right|<\delta } \varepsilon ^{-2/H}\Ex_\vartheta \left|{\bar\vartheta
  _\varepsilon-\vartheta} \right|^2\geq \Ex_{\vartheta_0} (\tilde{\xi}^2)= \Gamma ^{-\frac{2}{H}}\Ex (\tilde{\xi_o}^2). 
\end{equation}
Therefore we call the estimator $\vartheta _\varepsilon ^*$ asymptotically
efficient if for all $\vartheta _0\in \Theta $ we have the equality
\begin{equation}
\label{asef}
\lim_{\delta \rightarrow 0}\lim_{\varepsilon \rightarrow 0}
\sup_{\left|\vartheta -\vartheta _0\right|<\delta } \varepsilon ^{-2/H}\Ex_\vartheta \left|{\vartheta
  _\varepsilon^*-\vartheta} \right|^2= \Ex_{\vartheta_0} (\tilde{\xi}^2). 
\end{equation}
   The proof of this bound follows from the general results presented in
   \cite{IH81}. We can recall here the scetch of the proof supposing that the
   properties of the bayesian estimators for this model are already proved
   (see theorem \ref{T1} below). Introduce a continuous positive density
   function $\left(q\left(\vartheta  \right),\vartheta _0-\delta
   <\vartheta<\vartheta _0+\delta  \right) $. Then we can write
\begin{align*}
\sup_{\left|\vartheta -\vartheta _0\right|<\delta } \Ex_\vartheta \left|{\bar\vartheta
  _\varepsilon-\vartheta} \right|^2&\geq\int_{\vartheta
  _0-\delta}^{\vartheta _0+\delta}\Ex_\vartheta \left|{\bar\vartheta 
  _\varepsilon-\vartheta} \right|^2 q\left(\vartheta \right){\rm d}\vartheta\\
& \geq \int_{\vartheta
  _0-\delta}^{\vartheta _0+\delta}\Ex_\vartheta \left|{\tilde\vartheta 
  _{q,\varepsilon}-\vartheta} \right|^2 q\left(\vartheta \right){\rm d}\vartheta .
\end{align*}
where we denoted $\tilde\vartheta 
  _{q,\varepsilon} $ the bayesian estimator in the case of the density a
priory $q\left(\cdot \right)$.    As we have the convergence of moments of BE
we obtain  the limit
\begin{align*}
&\Liminf_{\varepsilon \rightarrow 0}
\sup_{\left|\vartheta -\vartheta _0\right|<\delta } \varepsilon ^{-2/H}\Ex_\vartheta \left|{\vartheta
  _\varepsilon^*-\vartheta} \right|^2\geq \lim_{\varepsilon \rightarrow 0}\varepsilon ^{-2/H}\int_{\vartheta
  _0-\delta}^{\vartheta _0+\delta}\Ex_\vartheta \left|{\tilde\vartheta 
  _{q,\varepsilon}-\vartheta} \right|^2 q\left(\vartheta \right){\rm
  d}\vartheta\\
&\qquad \qquad =\int_{\vartheta
  _0-\delta}^{\vartheta _0+\delta}\Ex_\vartheta | \tilde\xi  |^2 q\left(\vartheta \right){\rm
  d}\vartheta= \Ex ( \tilde\xi^2  )=\Gamma ^{-\frac{2}{H}}\Ex ( \tilde\xi_o  ^2)
\end{align*}
for all $\delta > 0$. Remind that $\Ex | \tilde\xi |^2$ does not depend on
$\vartheta $.  This proves the lower bound \eqref{lb}.

\begin{theorem}
\label{T1} The MLE and BE  are  consistent, have different limit distributions
\begin{align*}
\varepsilon ^{-\frac{1}{H}}\left(\hat\vartheta _\varepsilon-\vartheta
\right)\Longrightarrow \hat\xi ,\qquad \varepsilon ^{-\frac{1}{H}}\left(\tilde\vartheta
_\varepsilon-\vartheta 
\right)\Longrightarrow \tilde\xi ,
\end{align*}
 the polynomial moments converge : for any $p>0$
\begin{align*}
\lim_{\varepsilon \rightarrow 0} \Ex_\vartheta \left|\frac{\hat\vartheta
  _\varepsilon-\vartheta}{\varepsilon^{\frac{1}{H}} } \right|^p= \Ex_\vartheta
|\hat\xi |^p,\qquad \lim_{\varepsilon \rightarrow 0} \Ex_\vartheta
\left|\frac{\tilde\vartheta _\varepsilon-\vartheta}{\varepsilon ^{\frac{1}{H}}
} \right|^p= \Ex_\vartheta  |\tilde \xi |^p
\end{align*}
and the BE are asymptotically efficient. 
\end{theorem}
{\bf Proof.} To prove this theorem we check the conditions of the general
Theorem 1.10.1 in \cite{IH81}.  Let us put $\varphi _\varepsilon =\varepsilon
^{1/H}$ and introduce the normalized likelihood ratio
\begin{align*}
Z_\varepsilon \left(u\right)=\frac{V\left(\vartheta _0+\varphi _\varepsilon
  u,X^T\right)}{V\left(\vartheta _0,X^T\right)},\qquad u\in \UU_\varepsilon
=\left(\varepsilon ^{-1/H} \left(\alpha -\vartheta _0\right),\varepsilon
^{-1/H} \left(\beta -\vartheta _0\right)\right). 
\end{align*}
The verification of these conditions we do with the help of the lemmas
presented below. 
\begin{lemma}
\label{L0}
We have the
convergence of
finite-dimensional distributions of $ Z_\varepsilon \left(\cdot \right)$:
for any set $u_1,\ldots,u_k $ and any 
$k=1,2,\ldots$ 
\begin{equation}
\label{f-d}
\left( Z_\varepsilon \left(u_1\right),\ldots, Z_\varepsilon
\left(u_k\right)\right)\Longrightarrow  \left( Z \left(u_1\right),\ldots, Z
\left(u_k\right)\right).
\end{equation}
This convergence is uniforme in $\vartheta $ on compacts $\KK\subset \Theta
$. 
\end{lemma}
{\bf Proof.} 
 We can write ($u>0$)
\begin{align*}
\ln Z_\varepsilon \left(u\right)&=\frac{1}{\varepsilon
  ^2}\int_{0}^{T}\left[S\left(\vartheta _0+\varphi _\varepsilon
  u,t \right)-S\left(\vartheta _0,t \right)\right]{\rm d}X_t\\
& \qquad \qquad -\frac{1}{2\varepsilon
  ^2}\int_{0}^{T}\left[S\left(\vartheta _0+\varphi _\varepsilon
  u,t \right)^2-S\left(\vartheta _0,t \right)^2\right]{\rm d}t\\
&=\frac{1}{\varepsilon
  }\int_{0}^{T}\left[S\left(\vartheta _0+\varphi _\varepsilon
  u,t \right)-S\left(\vartheta _0,t \right)\right]{\rm d}W_t\\
&\qquad  \qquad -\frac{1}{2\varepsilon
  ^2}\int_{0}^{T}\left[S\left(\vartheta _0+\varphi _\varepsilon
  u,t \right)-S\left(\vartheta _0,t \right)\right]^2{\rm d}t.
\end{align*}
For the last integral we have
\begin{align*}
&\int_{0}^{T}\left[S\left(\vartheta _0+\varphi _\varepsilon
  u,t \right)-S\left(\vartheta _0,t \right)\right]^2{\rm d}t\\
&\quad =\int_{0}^{T}\left[a \left|t-\vartheta _0-\varphi _\varepsilon
  u\right|^\kappa  -a\left|t-\vartheta _0\right|^\kappa +h\left(\vartheta _0+\varphi _\varepsilon
  u,t \right)-h\left(\vartheta _0,t \right)\right]^2{\rm d}t\\
&\quad =\int_{-\vartheta _0}^{T-\vartheta _0}\left[a \left|t-\varphi _\varepsilon
  u\right|^\kappa  -a\left|t\right|^\kappa +\varphi _\varepsilon
  u\,\dot h(\tilde \vartheta,t-\vartheta _0 )\right]^2{\rm d}t,
\end{align*}
where we changed the variable and used Taylor expansion for the function
$h\left(\vartheta ,t\right)$.  

Let us put $t=\varphi _\varepsilon s$, then we obtain
\begin{align*}
&\int_{0}^{T}\left[S\left(\vartheta _0+\varphi _\varepsilon
  u,t \right)-S\left(\vartheta _0,t \right)\right]^2{\rm d}t\\
&\quad =\varphi _\varepsilon ^{2\kappa +1}\int_{-\frac{\vartheta _0}{\varphi
      _\varepsilon }}^{\frac{T-\vartheta _0}{\phi _\varepsilon }}\left[a
    \left|s-  u\right|^\kappa  -a\left|s\right|^\kappa +\varphi
    _\varepsilon^{1-\kappa } 
  u\,\dot h(\tilde \vartheta,s\varphi _\varepsilon -\vartheta _0 )\right]^2{\rm d}t\\
&\quad =a^2\varphi _\varepsilon ^{2\kappa +1}\int_{-\frac{\vartheta
      _0}{\varphi _\varepsilon }}^{\frac{T-\vartheta _0}{\phi _\varepsilon
  }}\left[ \left|s-  u\right|^\kappa  -\left|s\right|^\kappa \right]^2{\rm
    d}t\left(1+o\left(1\right)\right). 
\end{align*}
Hence
\begin{align}
&\frac{1}{\varepsilon ^2}\int_{0}^{T}\left[S\left(\vartheta _0+\varphi _\varepsilon
  u,t \right)-S\left(\vartheta _0,t \right)\right]^2{\rm d}t\nonumber\\
&\quad =\frac{a^2\varphi _\varepsilon ^{2\kappa +1}}{\varepsilon
    ^2}\int_{-\frac{\vartheta _0}{\varphi _\varepsilon }}^{\frac{T-\vartheta
      _0}{\varphi  _\varepsilon }}\left[ \left|s-  u\right|^\kappa
    -\left|s\right|^\kappa \right]^2{\rm d}s\left(1+o\left(1\right)\right)\nonumber\\ 
&\quad ={a^2\left|u\right|^{2\kappa +1}}\int_{-\frac{\vartheta _0}{\varphi
      _\varepsilon u
  }}^{\frac{T-\vartheta _0}{\varphi  _\varepsilon u}}\left[ \left|v-
    1\right|^\kappa  -\left|v\right|^\kappa \right]^2{\rm
    d}v\left(1+o\left(1\right)\right)\longrightarrow  \Gamma
  ^2\;\left|u\right|^{2\kappa +1}, 
\label{rep}
\end{align}
where we put $s=vu$. 

The similar  calculations for the  stochastic integral provide us the relations
\begin{align*}
&\frac{1}{\varepsilon }\int_{0}^{T}\left[S\left(\vartheta _0+\varphi _\varepsilon
  u,t \right)-S\left(\vartheta _0,t \right)\right]{\rm d}W_t \\
&\quad =a\int_{-\frac{\vartheta _0}{\varphi _\varepsilon}}^{\frac{T-\vartheta
      _0}{\varphi  _\varepsilon }}\left[ \left|s- 
    u\right|^\kappa  -\left|s\right|^\kappa \right]{\rm
    d}\tilde W\left(s\right)\left(1+o\left(1\right)\right)\\
&\quad \Longrightarrow a\int_{-\infty }^{\infty }\left[ \left|s-
    u\right|^\kappa  -\left|s\right|^\kappa \right]{\rm
    d} W\left(s\right)\sim {\cal N}\left(0, \left|u\right|^{2H}\Gamma^2 \right).
\end{align*}
Here $  W\left(v\right),u\in R $ is two-sided Wiener process
\begin{align*}
W\left(v\right)=\begin{cases}
W_+\left(v\right),& \text{if $v>0$},\\
W_-\left(-v\right),& \text{if $v\leq 0$},
\end{cases}
\end{align*}
where $ W_+\left(v\right),W_-\left(v\right),v\geq 0$ are two independent
Wiener processes.

Let us denote 
\begin{align*}
W^H\left(u\right)=\Gamma ^{-1}\int_{-\infty }^{\infty }\left[ \left|s-
    u\right|^\kappa  -\left|s\right|^\kappa \right]{\rm
    d} W\left(s\right)
\end{align*}
and verify \eqref{fBm}. We use below the equality
$ab=\frac{1}{2}\left[a^2+b^2-\left(a-b\right)^2\right]$ 
\begin{align*}
& \Ex W^H\left(u\right)W^H\left(v\right)\\
&\qquad =\frac{1}{2}\left[ \Ex
  \left(W^H\left(u\right)\right)^2 + \Ex
  \left(W^H\left(v\right)\right)^2 -\Ex
  \left(W^H\left(u\right)-  W^H\left(v\right)\right)^2 \right]\\
&\qquad =\frac{1}{2} \left[\left|u\right|^{2H}+\left|v\right|^{2H}   -\left|u-v\right|^{2H}\right]
\end{align*}
because
\begin{align*}
&\Ex  \left(W^H\left(u\right)-  W^H\left(v\right)\right)^2 =\Gamma
^{-2}\int_{-\infty }^{\infty }\left[ \left|s-u\right|^\kappa
  -\left|s-v\right|^\kappa \right]^2{\rm d}s\\
&\qquad \qquad =\Gamma
^{-2}\int_{-\infty }^{\infty }\left[ \left|r-\left(u-v\right)\right|^\kappa
  -\left|r\right|^\kappa \right]^2{\rm d}r =\left|u-v\right|^{2\kappa +1}.
\end{align*}
Hence $W^H\left(u\right),u\in R$ is a double-sided  fBm.

Therefore we proved the convergence of one-dimensional distributions. The
multi-dimensional case is treated by a similar way. We have to verify the
convergence 
\begin{align*}
\sum_{j=1}^{k}\lambda _j\ln Z_\varepsilon \left(u_j\right)\Longrightarrow
\sum_{j=1}^{k}\lambda _j\ln Z \left(u_j\right). 
\end{align*}
for an arbitrary vectors $(\lambda _1,\ldots,\lambda _k)$ and $(u _1,\ldots,u _k)$. 

Let us denote
\begin{align*}
\Phi \left(\vartheta,\vartheta _0 \right)=\int_{0}^{T}\left[S\left(\vartheta
  ,t\right)-S\left(\vartheta_0 ,t\right)\right]^2{\rm d}t . 
\end{align*}
We have the following elementary estimate
\begin{lemma}
\label{L1} There exists a constant $\mu   >0$ such that
\begin{equation}
\label{ng}
\Phi \left(\vartheta,\vartheta _0 \right)\geq \mu  \left|\vartheta
- \vartheta_0 \right|^{2H}. 
\end{equation}
\end{lemma}
{\bf Proof.} Note that  for any $\nu >0$
$$
m\left(\nu \right)=\inf_{\left|\vartheta -\vartheta_0 \right|>\nu }\Phi
\left(\vartheta,\vartheta _0 \right)>0. 
$$ 
Indeed, if  for some $\nu >0$ we have $m\left(\nu \right)=0 $, then there
exists $\vartheta _1\not=\vartheta _0$ such that for all $t\in \left[0,T\right]$
\begin{align*}
a\left|t-\vartheta _1\right|^\kappa +h\left(\vartheta_1 ,t\right)=a\left|t-\vartheta _0\right|^\kappa +h\left(\vartheta_0 ,t\right)
\end{align*}
and the function
\begin{align*}
h\left(\vartheta_1 ,t\right)=a\left|t-\vartheta _0\right|^\kappa-a\left|t-\vartheta _1\right|^\kappa  +h\left(\vartheta_0 ,t\right)
\end{align*}
has no continuous bounded derivativ on $\vartheta _1$. 
Hence for $\left|\vartheta -\vartheta _0\right|>\nu $
\begin{align*}
\Phi \left(\vartheta,\vartheta _0 \right)\geq m\left(\nu
\right)\geq m\left(\nu\right)\frac{  \left|\vartheta -\vartheta_0\right|
  ^{2H}}{\left|\beta -\alpha\right| ^{2H}} .
\end{align*}

Further, for the values $\left|\vartheta -\hat\vartheta \right|\leq \nu $ for
sufficiently small $\nu $  we have
\begin{align*}
\Phi \left(\vartheta,\vartheta _0 \right)=\left|\vartheta -\vartheta _0\right|^{2H} \Gamma ^2\left(1+o\left(1\right)\right).
\end{align*}
Therefore for sufficiently small $\nu $ we can write 
\begin{align*}
\Phi \left(\vartheta,\vartheta _0 \right)\geq \frac{1}{2}\Gamma^2 \left|\vartheta -\vartheta _0\right| ^{2H}. 
\end{align*}
Taking 
\begin{align*}
\mu  =\min \left(  \frac{m\left(\nu\right)}{\left|\beta -\alpha\right| ^{2H}}
,\frac{\Gamma^2 }{2}
\right)
\end{align*}
we obtain \eqref{ng}.

This estimate allows us to verify the boundness of all moments of the pseudo
likelihood ratio process.

\begin{lemma}
\label{L2}
 There exist a constant $c>0$  such that
\begin{equation}
\label{tails}
\Ex_{\vartheta _0}  Z_\varepsilon ^{\frac{1}{2}}\left(u\right)\leq e^{-c \left|u\right|^{2H}}.
\end{equation}
\end{lemma}
{\bf Proof.}  We have
\begin{align*}
\Ex_{\vartheta _0} Z_\varepsilon
\left(u\right)^{\frac{1}{2}}=\exp\left\{-\frac{1}{8\,\varepsilon
  ^{2}}\Phi(\vartheta_0+\varphi _\varepsilon u,\vartheta _0)\right\}\leq
\exp\left\{-\frac{\mu }{8} \left|u\right|^{2H}  \right\},
\end{align*}
where we used \eqref{ng}.

\begin{lemma}
\label{L3} For any $N>0$ and $\left|u_1\right|<N$, $\left|u_2\right|<N$ we have
the estimate
\begin{equation}
\label{rr}
\Ex_{\vartheta _0}\left|  Z_\varepsilon ^{\frac{1}{2}}\left({u_2}\right)-
Z_\varepsilon ^{\frac{1}{2}}\left({u_1}\right)\right|^2 \leq C\left(1+N\right)\left|u_2-u_1\right|^{2H}
\end{equation}
with some constant $C>0$.
\end{lemma}
{\bf Proof.} We can write 
\begin{align*}
&\Ex_{\vartheta _0}\left|  Z_\varepsilon ^{\frac{1}{2}}\left({u_2}\right)-
Z_\varepsilon ^{\frac{1}{2}}\left({u_1}\right)\right|^2
=2\left(1-\Ex_{\vartheta _0+\varphi _\varepsilon u_1}
\left(\frac{Z_\varepsilon \left({u_2}\right) }{Z_\varepsilon
  \left({u_1}\right)}\right)^{\frac{1}{2}}\right)\\ 
&\qquad \qquad =2\left(1-\exp\left\{-\frac{1}{8\varepsilon ^2}
\Phi\left(\vartheta _0+\varphi _\varepsilon u_2,\vartheta _0+\varphi
_\varepsilon u_1 \right)\right\}\right)\\
&\qquad \qquad \leq \frac{1}{4\varepsilon ^2}
\Phi\left(\vartheta _0+\varphi _\varepsilon u_2,\vartheta _0+\varphi
_\varepsilon u_1 \right)\\
&\qquad \qquad =\frac{1}{4\varepsilon ^2}\int_{0}^{T} \left[
  a\left|t-\vartheta _0-\varphi _\varepsilon u_2\right|^\kappa
-a\left|t-\vartheta _0-\varphi _\varepsilon u_1\right|^\kappa\right.\\
&\qquad \qquad\quad \qquad \qquad\left. +
  h\left(\vartheta _0+\varphi _\varepsilon u_2,t\right)-h\left(\vartheta
  _0+\varphi _\varepsilon u_1,t\right) \right]^2{\rm d}t \\
&\qquad \qquad\leq \frac{1}{2\varepsilon ^2}\int_{0}^{T} \left[
  a\left|t-\vartheta _0-\varphi _\varepsilon u_2\right|^\kappa
-a\left|t-\vartheta _0-\varphi _\varepsilon u_1\right|^\kappa\right]^2{\rm d}t\\
&\qquad \qquad\quad \qquad \qquad + \frac{1}{2\varepsilon ^2}\int_{0}^{T} \left[ 
  h\left(\vartheta _0+\varphi _\varepsilon u_2,t\right)-h\left(\vartheta
  _0+\varphi _\varepsilon u_1,t\right) \right]^2{\rm d}t\\
&\qquad \qquad\leq \frac{\varphi _\varepsilon ^{2\kappa +1}}{2\varepsilon
  ^2}\Gamma ^2\left|u_2-u_1\right| ^{2\kappa +1}+\frac{\varphi _\varepsilon ^{2}}{2\varepsilon
  ^2} \int_{0}^{T} \dot h(\tilde \vartheta ,t)^2{\rm d}t\left(u_2-u_1\right)^2\\
&\qquad \qquad\leq C\left(1+\left|u_2-u_1\right|^{1-2\kappa }\right)\left|u_2-u_1\right|
^{2\kappa +1}\leq C\left(1+N\right)\left|u_2-u_1\right|^{2H} .
\end{align*}

Note that $2\kappa <1$ and $2H>1$. 
The properties of the likelihood ratio \eqref{f-d}, \eqref{tails} and
\eqref{rr} correspond to the conditions  of the Theorems 1.10.1 and 1.10.2 in
\cite{IH81} and therefore the MLE $\hat\vartheta _\varepsilon $ and BE
$\tilde\vartheta _\varepsilon $ have  all 
mentioned in the Theorem \ref{T1} properties.

{\bf Remark 2.1.} More detailed analysis shows that if the signal has several
points of cusp, say
\begin{align*}
S\left(\vartheta ,t\right)=\sum_{l=1}^{L}a_l\left|t-\vartheta \right|^{\kappa
  _l}, 
\end{align*}
where $\kappa _l\in (0,\frac{1}{2})$, then the result of the Theorem \ref{T1} holds with 
\begin{align*}
\kappa =\min_{1\leq l\leq L} \kappa _l.
\end{align*}
The proof is similar to the given proof of the Theorem \ref{T1}.

{\bf Remark 2.2.} It is possible to study the properties of the estimators
$\hat\vartheta _\varepsilon  $ and $\tilde\vartheta _\varepsilon $ in the case
of multiple different singularities. For example, suppose that
\begin{align*}
S\left(\vartheta ,t\right)=\sum_{l=1}^{L}a_l\left|t-\vartheta_l \right|^{\kappa
  _l}, 
\end{align*}
where $\vartheta =\left(\vartheta _1,\ldots,\vartheta _L\right)\in \Theta
$. Here
$\Theta =\left(\alpha _1,\beta _1\right)\times \cdots \times \left(\alpha
_L,\beta _L\right)$, $0<\alpha _l<\beta _l<T$ and  $\beta _l<\alpha _{l+1}$, $ l=1,\ldots, L-1$. 

Then the limit for the normalized likelihood ratio 
\begin{align*}
Z_\varepsilon \left(u_1,\ldots,u_L\right)=\frac{V\left(\vartheta
  _l+\varepsilon ^{\frac{1}{H_1}}u_1,\ldots, \vartheta
  _L+\varepsilon ^{\frac{1}{H_L} }u_L,X^T\right)}{V\left(\vartheta
  _l,\ldots, \vartheta
  _L,X^T\right)}
\end{align*}
 is the process
\begin{align*}
Z\left(u_1,\ldots,u_L\right)=\prod_{l=1}^L Z_l\left(u_l\right),\qquad u_l\in R,
\end{align*}
where
\begin{align*}
%\label{z}
Z_l\left(u_l\right)=\exp \left\{\Gamma_l W_l^{H_l}\left(u_l\right)-\frac{\Gamma_l
   ^2}{2}\left|u_l\right|^ {2H_l}\right\},\qquad u_l\in R,
\end{align*}
and the constants
\begin{align*}
\Gamma_l^2={a_l^2}\int_{-\infty }^{\infty }\left[ \left|v-
    1\right|^{\kappa_l}  -\left|v\right|^{\kappa_l }\right]^2{\rm
    d}v.
\end{align*}
The fBm processes $\left(W_l^{H_l}\left(\cdot \right),\ldots,
W_L^{H_L}\left(\cdot \right) \right)$ are independent. The MLE $\hat\vartheta
_\varepsilon=\left(\hat\vartheta _{1,\varepsilon },\ldots, \hat\vartheta
_{L,\varepsilon } \right) $  and BE $\tilde\vartheta _\varepsilon =\left(\tilde\vartheta _{1,\varepsilon },\ldots, \tilde\vartheta
_{L,\varepsilon } \right)$ are
defined by the same relations \eqref{mle}, \eqref{be} and have different 
 limit distributions. In particularly, for the MLE we have the convergence 
\begin{align*}
\left(\frac{\hat\vartheta _{1,\varepsilon }-\vartheta _1}{\varepsilon
  ^{\frac{1}{H_1}}},\ldots,\frac{ \hat\vartheta _{L,\varepsilon }-\vartheta _L}{\varepsilon
  ^{\frac{1}{H_L}}}   \right)\Longrightarrow \left(\hat\xi _1,\ldots,\hat\xi _L \right).
\end{align*}
The limit random variables $ \left(\hat\xi _1,\ldots,\hat\xi _L \right)$ are
defined by the equations
\begin{align*}
Z_l(\hat\xi _l)=\sup_u Z_l\left(u\right),\qquad l=1,\ldots,L
\end{align*}
and are  independent. Of course, the
bayesian estimators have the same rate and the  asymptotic distribution is
\begin{align*}
\left(\frac{\tilde\vartheta _{1,\varepsilon }-\vartheta _1}{\varepsilon
  ^{\frac{1}{H_1}}},\ldots,\frac{ \tilde\vartheta _{L,\varepsilon }-\vartheta _L}{\varepsilon
  ^{\frac{1}{H_L}}}   \right)\Longrightarrow \left(\tilde\xi _1,\ldots,\tilde\xi _L \right).
\end{align*} 
Here the random variables 
\begin{align*}
\tilde\xi _l=\frac{\int_{}^{}u_lZ_l\left(u_l\right){\rm
    d}u_l}{\int_{}^{}Z_l\left(u_l\right){\rm d}u_l} ,\qquad l=1,\ldots,L
\end{align*}
are as well asymptotically independent.

\section{Misspecification}

We are intrerested by the following problem of misspecification. Suppose that the
model of observations  choosen by the statistician 
 ({\it theoretical model}) is 
\begin{align*}
{\rm d}X_t=M\left(\vartheta ,t\right){\rm d}t+\varepsilon {\rm d}W_t,\quad
X_0=0,\quad 0\leq t\leq T. 
\end{align*}
The signal $M\left(\vartheta ,t\right)$  is supposed to be 
$$
M\left(\vartheta ,t \right) =a\left|t-\vartheta \right|^\kappa
,\qquad 0\leq t\leq T, 
$$
where $\kappa \in (0,\frac{1}{2})$ and $\vartheta \in \Theta = \left(\alpha
<\vartheta <\beta \right)$. As before we suppose that $0<\alpha <\beta <T$. 

 The observed  process ({\it real model})  is 
\begin{align}
\label{11}
{\rm d}X_t=S\left(\vartheta_0 ,t\right){\rm d}t+\varepsilon {\rm d}W_t,\quad
X_0=0,\quad 0\leq t\leq T, 
\end{align}
where $\vartheta _0\in \Theta $ is the true value   and the function $S\left(\vartheta ,\cdot
\right)\in L_2\left(0,T\right)$ is sufficiently smooth.

  The 
likelihood ratio function (misspecified) is 
\begin{align*}
%\label{12}
V\left(\vartheta ,X^T\right)=\exp\left\{\frac{1}{\varepsilon ^2}\int_{0}^{T}
M\left(\vartheta ,t \right) {\rm d}X_t-\frac{1}{2\varepsilon ^2}\int_{0}^{T}
M\left(\vartheta ,t \right)^2 {\rm d}t \right\},\; \vartheta \in \Theta 
\end{align*}
where we have to substitute the observations from the equation
\eqref{11}. Therefore  the (pseudo) MLE $\hat\vartheta _\varepsilon $ is defined by the equation
\begin{align}
\label{ss}
V\left(\hat\vartheta _\varepsilon ,X^T\right)=\sup_{\vartheta \in \Theta
}V\left(\vartheta ,X^T\right). 
\end{align}

To see the limit of the MLE we write the likelihood ratio as
follows
\begin{align*}
&\varepsilon ^2\ln V\left(\vartheta ,X^T\right)\\
&\qquad=\varepsilon \int_{0}^{T}
M\left(\vartheta ,t \right) {\rm d}W_t -\frac{1}{2}\int_{0}^{T}
\left[M\left(\vartheta ,t \right)^2 -2M\left(\vartheta ,t
  \right)S\left(\vartheta _0,t\right)\right]{\rm d}t \\
&\qquad=\varepsilon \int_{0}^{T}
M\left(\vartheta ,t \right) {\rm d}W_t -\frac{1}{2}
\left\|M\left(\vartheta ,\cdot \right) -S\left(\vartheta _0,\cdot\right)
\right\|^2+\frac{1}{2}\left\|S\left(\vartheta _0,\cdot \right) \right\|^2 
\end{align*}
where we denoted as $ \left\|\cdot \right\|$ the $L_2\left(0,T\right)$ norm. 
It is easy  to verify  the convergence
\begin{align*}
\sup_{\vartheta \in \Theta }\left|\varepsilon ^2\ln L\left(\vartheta ,X^T\right)- \frac{1}{2}
\left\|M\left(\vartheta ,\cdot \right) -S\left(\vartheta _0,\cdot\right)
\right\|^2+\frac{1}{2}\left\|S\left(\vartheta _0,\cdot \right) \right\|^2  \right|\rightarrow 0.
\end{align*}
 Suppose that the equation 
\begin{align*}
\inf_{\vartheta} \left\|M(\vartheta ,\cdot ) -S\left(\vartheta _0,\cdot\right)
\right\|=\left\|M(\hat\vartheta ,\cdot ) -S\left(\vartheta _0,\cdot\right)
\right\|
\end{align*}
has a unique solution $\hat\vartheta \in \Theta $.

Then we obtain as usual in such situations  that
the MLE $\hat\vartheta _\varepsilon $
converges to the value $\hat\vartheta $, which minimizes the Kullback-Leibler
distance.

 It is interesting to note that in general case $\hat\vartheta
\not=\vartheta _0$ but sometimes $\hat\vartheta =\vartheta
_0$ and we consider the conditions of the consistency in such situations. The
most interesting for us is the question of the rate of convergence of the MLE
to the true value.

Introduce the function
\begin{align*}
\Phi (\vartheta ,\hat\vartheta )=\left\|M(\vartheta ,\cdot )
-S\left(\vartheta _0,\cdot\right) \right\|^2-\|M(\hat\vartheta ,\cdot ) -S\left(\vartheta _0,\cdot\right)
\|^2
\end{align*}
and the conditions of regularity:

{\it Condition ${\cal M}$}.
\begin{enumerate}
\item {\it The parameter $\kappa \in \left(0,\frac{1}{2}\right)$.
\item The function $S\left(\vartheta _0,t\right)$ for all $\vartheta _0\in
  \Theta $ is two times continuously 
  differentiable w.r.t. $t\in \left[0,T\right]$.
\item The function $\Phi (\vartheta ,\hat\vartheta ) $ for all $\vartheta _0\in
  \Theta $  has a unique   minimum at the point
  $\hat\vartheta=\hat\vartheta\left(\vartheta _0\right)$. 
\item It's second derivative 
\begin{align*}
\gamma(\hat\vartheta )\equiv \left.\frac{\partial^2 \Phi(\vartheta
  ,\hat\vartheta ) }{\partial \vartheta ^2}\right|_{\vartheta
  =\hat\vartheta } >0
\end{align*}
 for all $\vartheta _0\in   \Theta $. }

\end{enumerate}

Let us denote 
\begin{align*}
&\hat Z\left(u\right)=\exp\left\{aW^H\left(u\right)-\frac{\gamma(\hat\vartheta
  )}{4} u^{2}\right\},\quad u\in R\\
& \hat
Z^o\left(u\right)=\exp\left\{w^H\left(v\right)-\frac{ v^{2}}{2}\right\},\qquad 
v\in R
\end{align*}
and define the random variables $\hat \zeta,\hat \zeta_o  $ by the relations
\begin{align*}
\hat Z(\hat \zeta  )=\sup_u \hat  Z\left(u\right),\qquad \quad \hat Z^o(\hat \zeta_o
)=\sup_v \hat  Z^o\left(v\right). 
\end{align*}
Note that 
\begin{equation}
\label{eqi}
\hat \zeta =\left(\frac{2a}{\gamma (\hat\vartheta )}\right)^{\frac{H}{2H-1}}\;\hat \zeta_o.
\end{equation}
To verify \eqref{eqi} we change the variables $u=rv$ with
$r=\left(2a\right)^{\frac{H}{2H-1}}\gamma (\hat\vartheta )^{-\frac{H}{2H-1}}
$ and write
\begin{align*}
&aW^H\left(u\right)-\frac{\gamma(\hat\vartheta
  )}{4} u^{2}=aW^H\left(rv\right)-\frac{\gamma(\hat\vartheta
  )r^2}{4} v^{2}\\
&\qquad \quad =ar^{\frac{1}{H}}\left(\frac{W^H\left(rv\right) }{
  r^{\frac{1}{H}} }-\frac{\gamma(\hat\vartheta 
  )r^{2-\frac{1}{H}}}{4a} v^{2}
  \right)=ar^{\frac{1}{H}}\left(w^H\left(v\right)-\frac{v^2}{2}\right), 
\end{align*}
where the fBm $w^H\left(v\right)= r^{-\frac{1}{H}} W^H\left(rv\right)$.
\begin{theorem}
\label{T2} Let the conditions ${\cal M}$ be fulfilled, then the estimator
$\hat\vartheta _\varepsilon $  converges to the value $\hat\vartheta $,  
has the limit distribution  
\begin{align}
\label{asno}
\frac{\hat\vartheta _\varepsilon-\hat\vartheta }{\varepsilon
  ^{\frac{2}{3-2\kappa }}}\Longrightarrow \hat \zeta ,
\end{align}
and  for any $p>0$
\begin{align}
\label{conv}
\lim_{\varepsilon \rightarrow 0}\Ex_{\vartheta } \left|\frac{\hat\vartheta
  _\varepsilon-\hat\vartheta }{\varepsilon ^{\frac{2}{3-2\kappa }}}\right|^p = \Ex_{\vartheta
  }\left|\hat \zeta\right|^p=\left(\frac{2a}{\gamma (\hat\vartheta
  )}\right)^{\frac{pH}{2H-1}}\Ex\left|\hat \zeta_o\right|^p.
\end{align}

\end{theorem}
{\bf Proof.} Introduce the normalized pseudo-likelihood ratio process
\begin{align*}
Z_\varepsilon \left(u\right)=\frac{V\left(\hat\vartheta +\varphi _\varepsilon
  u,X^T\right)}{V\left(\hat\vartheta,X^T \right)} ,\qquad
u\in \UU_\varepsilon =\left(   \frac{\left(\alpha -\hat\vartheta
  \right)}{\varphi _\varepsilon },\frac{\left(\beta  -\hat\vartheta
  \right)}{\varphi _\varepsilon }\right),  
\end{align*}
where $\varphi _\varepsilon\rightarrow 0$ will be defined later and denote $\vartheta
_u=\hat\vartheta+\varphi _\varepsilon   u  $. Below we use the same arguments
as that of  the preceding section in similar situation ($u>0$)
\begin{align*}
&\ln Z_\varepsilon \left(u\right)=\frac{1}{\varepsilon ^2}\int_{0}^{T} \left[
  M\left(\hat\vartheta +\varphi _\varepsilon u,t\right)-M\left(\hat\vartheta
  ,t\right)\right]{\rm d}X_t\\
&\qquad \quad \qquad -\frac{1}{2\varepsilon ^2}\int_{0}^{T}  \left[
  M\left(\hat\vartheta +\varphi _\varepsilon u,t\right)^2-M\left(\hat\vartheta
  ,t\right)^2\right]{\rm d}t\\
&\qquad=\frac{1}{\varepsilon
}\int_{0}^{T}\left[a \left|t- \vartheta_u \right|^\kappa -a\left|t-
    \hat\vartheta \right|^\kappa\right] {\rm d}W_t 
  \\ 
&\qquad \quad\qquad -\frac{1}{2\varepsilon^2
}\int_{0}^{T}\left[ a\left|t- \vartheta_ u \right|^\kappa -a\left|t-
    \hat\vartheta \right|^\kappa+ h(\vartheta_u,t)-h(\hat\vartheta
    ,t)\right]\\
&\qquad \quad\qquad \quad 
 \left[a \left|t- \vartheta_ u \right|^\kappa +a\left|t- \hat\vartheta
    \right|^\kappa- 2S(\vartheta_0,t)\right] {\rm d}t\\
&\qquad=\frac{a \varphi _\varepsilon ^{\kappa +\frac{1}{2}}}{\varepsilon
}\int_{-\frac{\hat \vartheta }{\varphi _\varepsilon}}^{\frac{T-\hat \vartheta
   }{\varphi _\varepsilon}}\left[ \left|s- u \right|^\kappa
   -\left|s\right|^\kappa\right] {\rm d}W\left(s\right) 
  -\frac{1}{2\varepsilon^2}\;\Phi(\vartheta _u,\hat\vartheta  ) .
\end{align*}

Let us study the function $\Phi (\vartheta _u,\hat\vartheta
  ) $ for a fixed $u>0$ as $\varphi _\varepsilon \rightarrow 0$. We have
\begin{align*}
\Phi (\vartheta ,\hat\vartheta  )&=\int_{0}^{T} \left[
  M\left(\vartheta ,t\right)-S\left(\vartheta _0,t\right)\right] ^2{\rm d}t-\int_{0}^{T} \left[
  M(\hat\vartheta,t)-S\left(\vartheta _0,t\right)\right] ^2{\rm
  d}t\\
&=\int_{-\vartheta}^{T-\vartheta} \left[a\left|s\right|^\kappa
  -S\left(\vartheta _0,s+\vartheta \right) 
  \right] ^2{\rm d}t-\int_{0}^{T} \left[
  M(\hat\vartheta,t)-S\left(\vartheta _0,t\right)\right] ^2{\rm
  d}t
\end{align*}
and
\begin{align*}
\Phi'_\vartheta  (\vartheta ,\hat\vartheta  )&=\left[a\left| \vartheta  \right|^\kappa
  -S\left(\vartheta _0,0 \right) 
  \right] ^2  -\left[a\left| T-\vartheta  \right|^\kappa
  -S\left(\vartheta _0,T \right) 
  \right] ^2 \\
&\qquad -2  \int_{-\vartheta }^{T-\vartheta } \left[
  a\left|s\right|^\kappa -S\left(\vartheta _0,s+\vartheta\right)\right]S'\left(\vartheta
_0,s+\vartheta\right) {\rm d}s .
\end{align*}
 Recall that as $\hat\vartheta\in \Theta  $ is
the point of minimum of the function $\Phi  (\vartheta ,\hat\vartheta
),\vartheta \in \Theta  $ we have the equalities
\begin{align*}
\Phi   (\hat\vartheta ,\hat\vartheta)=0 ,\qquad \quad 
\Phi'_\vartheta   (\hat\vartheta ,\hat\vartheta
)=0 .
\end{align*}
Let us write the Taylor expansion
\begin{align*}
\Phi   (\vartheta_u ,\hat\vartheta)&=\Phi   (\hat\vartheta
,\hat\vartheta)+\varphi _\varepsilon u \Phi'_\vartheta   (\hat\vartheta
,\hat\vartheta)+\frac{\varphi _\varepsilon^2 u^2}{2} 
\Phi''_\vartheta   (\hat\vartheta
,\hat\vartheta)\left(1+o\left(1\right)\right)\\
&=\frac{\varphi _\varepsilon^2 u^2}{2} 
\Phi''_\vartheta   (\hat\vartheta ,\hat\vartheta)\left(1+o\left(1\right)\right)
\end{align*}
and study the difference
\begin{align*}
&\Phi'_\vartheta   (\vartheta_u,\hat\vartheta)-\Phi'_\vartheta   (\hat\vartheta
,\hat\vartheta)=\left[a\left| \vartheta_u  \right|^\kappa 
  -S\left(\vartheta _0,0 \right) 
  \right] ^2 -\left[a| \hat\vartheta  |^\kappa 
  -S\left(\vartheta _0,0 \right) 
  \right] ^2\\
 &\quad \quad  + \left[a| T-\hat\vartheta |^\kappa
  -S\left(\vartheta _0,T \right) 
  \right] ^2 -\left[a\left| T-\vartheta_u  \right|^\kappa
  -S\left(\vartheta _0,T \right) 
  \right] ^2 \\
&\qquad -2  \int_{-\vartheta_u }^{T-\vartheta_u } \left[
  a\left|s\right|^\kappa -S\left(\vartheta _0,s+\vartheta_u\right)\right]S'\left(\vartheta
_0,s+\vartheta_u\right) {\rm d}s\\
&\qquad +2  \int_{-\hat\vartheta }^{T-\hat\vartheta} \left[
  a\left|s\right|^\kappa -S\left(\vartheta _0,s+\hat\vartheta\right)\right]S'\left(\vartheta
_0,s+\hat\vartheta\right) {\rm d}s.
\end{align*}
We have the estimates
\begin{align*}
&\left[a\left| \vartheta_u \right|^\kappa -S\left(\vartheta _0,0 \right)
  \right] ^2 -\left[a| \hat\vartheta |^\kappa -S\left(\vartheta _0,0 \right)
  \right] ^2\\
&\qquad \quad= a\left[| \hat\vartheta +\varphi _\varepsilon u|^\kappa -| \hat\vartheta
    |^\kappa\right]  \left[a| \hat\vartheta +\varphi _\varepsilon u|^\kappa +a| \hat\vartheta
    |^\kappa -2 S\left(\vartheta _0,0 \right)
  \right] \\
&\qquad \quad= a\left[| \hat\vartheta +\varphi _\varepsilon u|^\kappa -| \hat\vartheta
    |^\kappa\right]  \left[a| \hat\vartheta +\varphi _\varepsilon u|^\kappa +a| \hat\vartheta
    |^\kappa -2 S\left(\vartheta _0,0 \right)
  \right] \\
&\qquad =\frac{2a\kappa }{\hat\vartheta^{1-\kappa } } \left[a| \hat\vartheta
    |^\kappa- S\left(\vartheta _0,0 \right) \right]      \varphi _\varepsilon
  u +O\left(\varphi _\varepsilon^2
  u^2\right)
\end{align*}
and similary
\begin{align*}
&\left[a\left|T- \hat\vartheta \right|^\kappa -S\left(\vartheta _0,T \right)
  \right] ^2 -\left[a| T-\hat\vartheta -\varphi _\varepsilon u|^\kappa -S\left(\vartheta _0,T \right)
  \right] ^2\\
&\qquad \quad =\frac{2a\kappa }{|T-\hat\vartheta |^{1-\kappa }}\left[a| T-\hat\vartheta |^\kappa-S\left(\vartheta _0,T \right) \right]\varphi _\varepsilon u+O\left(\varphi _\varepsilon^2
  u^2\right)
\end{align*}
because \begin{align*}
| \hat\vartheta +\varphi _\varepsilon u|^\kappa -| \hat\vartheta    |^\kappa&=|
\hat\vartheta    |^\kappa\left(1+ \frac{\kappa \varphi _\varepsilon
  u}{\hat\vartheta }\right)-| \hat\vartheta    |^\kappa +O\left(\varphi _\varepsilon^2
  u^2\right)\\
&=\frac{\kappa \varphi _\varepsilon
  u}{\hat\vartheta^{1-\kappa } }+O\left(\varphi _\varepsilon^2
  u^2\right).
\end{align*}
Further, we can write
\begin{align*}
& \int_{-\hat\vartheta }^{T-\hat\vartheta} \left[
  a\left|s\right|^\kappa -S\left(\vartheta _0,s+\hat\vartheta\right)\right]S'\left(\vartheta
_0,s+\hat\vartheta\right) {\rm d}s\\
&\qquad 
-2  \int_{-\vartheta_u }^{T-\vartheta_u } \left[
  a\left|s\right|^\kappa -S\left(\vartheta _0,s+\vartheta_u\right)\right]S'\left(\vartheta
_0,s+\vartheta_u\right) {\rm d}s \\
&\quad =\int_{0}^{T}\left|t-\hat\vartheta \right|^\kappa
\left[S'\left(\vartheta _0,     t+\varphi _\varepsilon u\right)-S'\left(\vartheta _0,t\right)
   \right] {\rm d}t \\
&\quad \quad +\int_{0}^{T}\left[ S\left(\vartheta
  _0,t+\varphi _\varepsilon u\right)S'\left(t+\varphi _\varepsilon u\right) -S\left(\vartheta _0,t\right)S'\left(\vartheta
_0,t\right) \right] {\rm d}t\\
&\quad \quad +\left(\int_{-\vartheta _u}^{-\hat\vartheta } -\int_{T-\vartheta _u}^{T-\hat\vartheta }\right)\left[
  a\left|s\right|^\kappa -S\left(\vartheta _0,s+\vartheta_u\right)\right]S'\left(\vartheta
_0,s+\vartheta_u\right) {\rm d}s .
\end{align*}
Therefore we obtain the relations
\begin{align*}
&\int_{0}^{T}\left|t-\hat\vartheta \right|^\kappa \left[S'\left(\vartheta _0,
    t+\varphi _\varepsilon u\right)-S'\left(\vartheta _0,t\right) \right] {\rm
    d}t \\ 
&\qquad = \int_{0}^{T}\left|t-\hat\vartheta \right|^\kappa
  S''\left(\vartheta _0, t\right) {\rm d}t\; \varphi _\varepsilon u
  +O\left(\varphi _\varepsilon^2 u^2\right),\\ 
&\int_{0}^{T}\left[
    S\left(\vartheta _0,t+\varphi _\varepsilon u\right)S'\left(\vartheta
    _0,t+\varphi _\varepsilon u\right) -S\left(\vartheta
    _0,t\right)S'\left(\vartheta _0,t\right) \right] {\rm d}t\\ 
&\qquad =
  \frac{1}{2}\int_{0}^{T}\left[S\left(\vartheta _0,t\right)^2 \right]''_t\;{\rm
    d}t\;\varphi _\varepsilon u +O\left(\varphi _\varepsilon^2
  u^2\right),\\ 
&\int_{-\vartheta _u}^{-\hat\vartheta } \left[
    a\left|s\right|^\kappa -S\left(\vartheta
    _0,s+\vartheta_u\right)\right]S'\left(\vartheta _0,s+\vartheta_u\right)
  {\rm d}s\\ &\qquad=
  \left[a\left|\hat\vartheta\right|^\kappa-S\left(\vartheta
    _0,0\right)\right]S'\left(\vartheta _0,0\right) \; \varphi _\varepsilon u
  +O\left(\varphi _\varepsilon^2 u^2\right),\\
&\int_{T-\vartheta _u}^{T-\hat\vartheta }\left[
  a\left|s\right|^\kappa -S\left(\vartheta _0,s+\vartheta_u\right)\right]S'\left(\vartheta
_0,s+\vartheta_u\right) {\rm d}s\\ 
&\qquad=  \left[a\left|T-\hat\vartheta\right|^\kappa-S\left(\vartheta
    _0,T\right)\right]S'\left(\vartheta _0,T\right) \; \varphi _\varepsilon u
  +O\left(\varphi _\varepsilon^2 u^2\right).
\end{align*}
All these together allows us to write
\begin{align}
&\frac{\Phi'_\vartheta   (\vartheta_u,\hat\vartheta)}{\varphi _\varepsilon u}=\frac{2a\kappa }{\hat\vartheta^{1-\kappa } } \left[a| \hat\vartheta
    |^\kappa- S\left(\vartheta _0,0 \right) \right] +\frac{2a\kappa }{|T-\hat\vartheta |^{1-\kappa }}\left[a| T-\hat\vartheta |^\kappa-S\left(\vartheta _0,T \right) \right]\nonumber\\
&\qquad+\left[a\left|\hat\vartheta\right|^\kappa-S\left(\vartheta
    _0,0\right)\right]S'\left(\vartheta _0,0\right)+\left[a\left|T-\hat\vartheta\right|^\kappa-S\left(\vartheta
    _0,T\right)\right]S'\left(\vartheta _0,T\right) \nonumber\\
&\qquad + 2\int_{0}^{T}\left|t-\hat\vartheta \right|^\kappa
  S''\left(\vartheta _0, t\right) {\rm d}t+\int_{0}^{T}\left[S\left(\vartheta
    _0,t\right)^2 \right]''_t\;{\rm    d}t+O\left(\varphi _\varepsilon u
  \right).
\label{sd0}
\end{align}
Hence we obtain the following expression for second derivative
\begin{align}
\label{sd}
&\Phi''_\vartheta   (\hat\vartheta,\hat\vartheta)=\lim_{\varphi _\varepsilon
  \rightarrow 0}\frac{\Phi'_\vartheta
  (\vartheta_u,\hat\vartheta)-\Phi'_\vartheta
  (\hat\vartheta,\hat\vartheta)}{\varphi _\varepsilon u} \nonumber\\
&\quad =\frac{2a\kappa }{\hat\vartheta^{1-\kappa } } \left[a| \hat\vartheta
    |^\kappa- S\left(\vartheta _0,0 \right) \right] +\frac{2a\kappa }{|T-\hat\vartheta |^{1-\kappa }}\left[a| T-\hat\vartheta |^\kappa-S\left(\vartheta _0,T \right) \right]\nonumber\\
&\qquad+\left[a\left|\hat\vartheta\right|^\kappa-S\left(\vartheta
    _0,0\right)\right]S'\left(\vartheta _0,0\right)+\left[a\left|T-\hat\vartheta\right|^\kappa-S\left(\vartheta
    _0,T\right)\right]S'\left(\vartheta _0,T\right) \nonumber\\
&\qquad + 2\int_{0}^{T}\left|t-\hat\vartheta \right|^\kappa
  S''\left(\vartheta _0, t\right) {\rm d}t+\int_{0}^{T}\left[S\left(\vartheta
    _0,t\right)^2 \right]''_t\;{\rm    d}t.
\end{align}

Now the log-likelihood ratio has the representation
\begin{align*}
\ln Z_\varepsilon \left(u\right)&=\frac{a\varphi _\varepsilon ^{\kappa
    +\frac{1}{2}}}{\varepsilon }W^H\left(u\right)
\left(1+o\left(1\right)\right) -\frac{\varphi _\varepsilon ^2u^2}{4\varepsilon
  ^2}\Phi''_\vartheta (\hat\vartheta,\hat\vartheta)
\left(1+o\left(1\right)\right)\\
&=\frac{\varphi _\varepsilon ^{\kappa
    +\frac{1}{2}}}{\varepsilon }\left(a  W^H\left(u\right)
\left(1+o\left(1\right)\right) -\frac{\varphi _\varepsilon ^{\frac{3}{2}-\kappa }}{\varepsilon
  }    \Phi''_\vartheta (\hat\vartheta,\hat\vartheta) \frac{u^2}{4}
\left(1+o\left(1\right)\right)  \right).
\end{align*}
Therefore if we put
\begin{align*}
\frac{\varphi _\varepsilon ^{\frac{3}{2}-\kappa }}{\varepsilon  }  =1,\qquad
\varphi _\varepsilon =\varepsilon  ^{\frac{2}{3-2\kappa }},\qquad \hat Z_\varepsilon \left(u\right)=Z_\varepsilon \left(u\right)^{\varepsilon ^\frac{4\kappa -2}{3-2\kappa }},
\end{align*}
then we obtain the convergence of finite-dimensional distributions 
\begin{align*}
\left(\hat Z_\varepsilon \left(u_1\right),\ldots,\hat Z_\varepsilon \left(u_k\right)\right)\Longrightarrow \left(\hat Z\left(u_1\right),\ldots,\hat Z \left(u_k\right)\right)
\end{align*}
for any $k=1,2,\ldots$.

Using the same arguments as in the proofs of the lemmae \ref{L1}-\ref{L3} we
obtain the relations
\begin{align*}
&\Phi \left(\vartheta ,\hat\vartheta \right)\geq \mu \left(\vartheta
-\hat\vartheta \right)^2,\\
&\Ex_{\vartheta _0}\hat Z_\varepsilon^{\frac{1}{2}} \left(u\right)\leq
e^{-cu^2},\\
&\Ex_{\vartheta _0}\left[\hat Z_\varepsilon^{\frac{1}{2}}
  \left(u_2\right)-\hat Z_\varepsilon^{\frac{1}{2}}
  \left(u_1\right)\right]^2\leq C\left(1+N\right)\left| u_2-u_1\right|^2 
\end{align*}
Therefore once more the asymptotic properties of the pseudo-MLE $\hat\vartheta
_\varepsilon $ follow from the general result by Ibragimov and Kasminskii
\cite{IH81}, Theorem 1.10.1. 

Let us remind how the properties of  $\hat\vartheta
_\varepsilon $ are related with the convergence of the stochastic processes
$\hat Z_\varepsilon \left(\cdot \right)\Longrightarrow \hat Z\left(\cdot
\right)$:
we can write 
\begin{align}
\label{an}
&\Pb_{\vartheta _0}\left(\frac{\hat\vartheta _\varepsilon-\hat\vartheta
}{\varphi _\varepsilon }<x\right)=\Pb_{\vartheta _0}\left(\hat\vartheta
_\varepsilon< \hat\vartheta+\varphi _\varepsilon x\right) \nonumber\\ &\qquad
\qquad =\Pb_{\vartheta _0}\left\{ \sup_{\vartheta <\hat\vartheta+\varphi
  _\varepsilon x} V\left(\vartheta ,X^T\right)>\sup_{\vartheta \geq
  \hat\vartheta+\varphi _\varepsilon x} V\left(\vartheta ,X^T\right) \right\}
\nonumber\\ &\qquad \qquad =\Pb_{\vartheta _0}\left\{ \sup_{\vartheta
  <\hat\vartheta+\varphi _\varepsilon x} \frac{V\left(\vartheta ,X^T\right)}{
  V\left(\hat\vartheta ,X^T \right)}>\sup_{\vartheta \geq \hat\vartheta+\varphi
  _\varepsilon x} \frac{V\left(\vartheta ,X^T \right) }{V\left(\hat\vartheta
  ,X^T \right)} \right\}\nonumber \\ &\qquad \qquad =\Pb_{\vartheta _0}\left\{
\sup_{u <x,u\in\UU_\varepsilon } Z_\varepsilon \left(u\right)>\sup_{u \geq
  x,u\in\UU_\varepsilon } Z_\varepsilon \left(u\right)
\right\}\nonumber\\ &\qquad \qquad =\Pb_{\vartheta _0}\left\{ \sup_{u
  <x,u\in\UU_\varepsilon } \hat Z_\varepsilon \left(u\right)>\sup_{u \geq
  x,u\in\UU_\varepsilon } \hat Z_\varepsilon \left(u\right) \right\}
=\Pb_{\vartheta _0}\left(\hat u_\varepsilon <x\right) ,
\end{align}
where $\hat u_\varepsilon=\frac{\hat\vartheta _\varepsilon-\hat\vartheta
}{\varphi _\varepsilon } $ is defined by the relation
\begin{align*}
 \hat Z_\varepsilon \left(\hat u_\varepsilon\right)= \sup_{u\in
   \UU_\varepsilon }\hat Z_\varepsilon \left(u\right) .
\end{align*} 
Now from the convergence $\hat Z_\varepsilon \left(\cdot \right)\Rightarrow
\hat Z \left(\cdot \right)$ we obtain
\begin{align*}
&\Pb_{\vartheta _0}\left\{ \sup_{u
  <x,u\in\UU_\varepsilon } \hat Z_\varepsilon \left(u\right)>\sup_{u \geq
  x,u\in\UU_\varepsilon } \hat Z_\varepsilon \left(u\right)
\right\}\\
&\qquad \qquad \qquad \longrightarrow \Pb_{\vartheta _0}\left\{ \sup_{u 
  <x } \hat Z \left(u\right)>\sup_{u \geq
  x } \hat Z \left(u\right) \right\}=\Pb_{\vartheta _0}\left(\hat\zeta <x\right)
\end{align*}
(see the details in \cite{IH81}, Theorem 1.10.1).

{\bf Remark 3.1.} Of course, it is possible to considere slightly more general
problem with the signal 
\begin{align*}
S\left(\vartheta ,t\right)=a\left|t-\vartheta \right|^\kappa \1_{\left\{t<\vartheta\right\}
}+b\left|t-\vartheta \right|^\kappa \1_{\left\{t\geq \vartheta\right\} }+h\left(\vartheta
,t\right), 
\end{align*}
where $b>0$ and $h\left(\vartheta ,t\right)$ is some smooth function of $\vartheta $ and
$t$. As usual in singular estimation problems the limit likelihood ratio $
Z\left(\cdot \right)$ does
not depend on the function $h\left(\cdot ,\cdot \right)$ and the properties of
the pseudo-MLE are quite close to that of the presented in the Theorem \ref{T2}.

There are another interesting problems of misspecification {\it cusp vs
  discontinuous} and {\it  discontinuous vs  cusp}, which can be illustrated
by the following example. Suppose that we have two signals
\begin{align*}
S\left(\vartheta ,t\right)=-a\left|t-\vartheta \right|^\kappa \1_{\left\{t<\vartheta\right\}
}+a\left|t-\vartheta \right|^\kappa \1_{\left\{t\geq \vartheta\right\} }, 
\end{align*}
where $\kappa \in (0,\frac{1}{2})$ and
\begin{align*}
M\left(\vartheta ,t\right)=a\;\sgn \left(t-\vartheta \right).
\end{align*}
One problem is the estimation of the parameter $\vartheta $  in the situation,
where $S\left(\vartheta_0 ,t\right)$ is the observed  signal and
$M\left(\vartheta ,t\right) $ is supposed (theoretical)  signal. The second problem
corresponds to the situation where the
observed signal is $M\left(\vartheta_0 ,t\right) $ and the theoretical signal
is $S\left(\vartheta ,t\right)$. The both problems are studied in the
forthcomming paper.

\section{Estimation of the parameter $\kappa $}

Let us consider the problem of estimation of the parameter $\kappa\in
\left(k,K\right), 0< k<K<\infty   $ by
observations
\begin{align*}
{\rm d}X_t=a\left|t-\rho \right|^{\kappa_0} {\rm d}t+\varepsilon {\rm d}W_t,\quad
X_0=0,\quad 0\leq t\leq T, 
\end{align*}
where $a>0$ and $\rho\in \left(0,T\right)$ are some known parameters. The
likelihood-ratio function is
\begin{align*}
V\left(\kappa ,X^T\right)=\exp\left\{\int_{0}^{T}\frac{a
  \left|t-\rho\right|^{\kappa }}{\varepsilon ^2} {\rm d}X_t- \int_{0}^{T}\frac{a^2
  \left|t-\rho\right|^{2\kappa }}{2\varepsilon ^2} {\rm d} t\right\},\qquad  \kappa\in
\left(k,K\right)
\end{align*}
and the MLE $\hat\kappa _\varepsilon $ is solution of the equation
\begin{align*}
V\left(\hat\kappa _\varepsilon ,X^T\right)=\sup_{\kappa \in \left(k,K\right)}V\left(\kappa ,X^T\right)
\end{align*}
This is regular prolem with the Fisher information 
\begin{align*}
\II\left(\kappa \right)=a^2\int_{0}^{T} \left|t-\rho \right|^{2\kappa }\left(\ln
\left|t-\rho \right|\right)^2 {\rm d}t>0.
\end{align*}

It is easy to see that the identification condition 
\begin{align*}
\inf_{\left|\kappa -\kappa _0\right|>\nu } \int_{0}^{T} \left(\left|t-\rho\right|^{2\kappa }-
\left|t-\rho\right|^{\kappa _0}\right)^2 {\rm d}t>0.
\end{align*}
is fulfilled for any $\kappa _0$ and any $\nu >0$. 

Therefore the asymptotic normality 
\begin{align*}
\frac{\hat\kappa _\varepsilon -\kappa_0 }{\varepsilon }\Longrightarrow {\cal
  N}\left(0, \II\left(\kappa_0 \right)^{-1}\right)
\end{align*}
follows from the general theorem devoted to the parameter estimation in
regular families (see Theorem 3.1.1 in \cite{IH81}). Just note that the
normalized likelihood ratio
\begin{align*}
Z_\varepsilon^* \left(v\right)=\frac{V\left(\kappa _0+\varepsilon
  v,X^T\right)}{V\left(\kappa _0,X^T\right)} ,\qquad v\in \VV_\varepsilon
=\left(\frac{k-\kappa _0}{\varepsilon  },\frac{K-\kappa _0}{\varepsilon  }\right)
\end{align*}
converges to the process
\begin{align}
\label{z*}
Z^*\left(v\right)=\exp\left\{v\Delta -\frac{u^2}{2}\II\left(\kappa_0 \right)
\right\},\qquad v\in R,
\end{align}
where $\Delta \sim {\cal
  N}\left(0, \II\left(\kappa_0 \right)\right)$

It is interesting as well to consider the problem of two-dimensional parameter
$\vartheta =\left( \rho,\kappa \right)$ estimation. The likelihood-ratio
function is
\begin{align*}
V\left(\rho ,\kappa ,X^T\right)=\exp\left\{\int_{0}^{T}\frac{a
  \left|t-\rho\right|^{\kappa }}{\varepsilon ^2} {\rm d}X_t- \int_{0}^{T}\frac{a^2
  \left|t-\rho\right|^{2\kappa }}{2\varepsilon ^2} {\rm d} t\right\},\qquad
\vartheta \in \Theta,
\end{align*}
where $ \Theta =\left(\alpha ,\beta \right)\times \left(k,K\right)$, $0<\alpha
<\beta <T$. 

It can be shown that the normalized likelihood ratio
\begin{align*}
Z_\varepsilon \left(u,v\right)=\frac{V\left(\rho _0+\varepsilon ^{\frac{1}{H}}u,\kappa _0+\varepsilon
  v,X^T\right)}{V\left(\rho _0,\kappa _0,X^T\right)} 
\end{align*} 
converges to the random process
\begin{align*}
Z\left(u,v\right)=Z\left(u\right)\,Z^*\left(v\right)
\end{align*}
where the processes $Z\left(\cdot \right) $ and   $Z^*\left(\cdot \right) $
are defined by the expressions \eqref{z} and \eqref{z*}. Note that the fBm $W^H\left(\cdot \right)$
 and the random variable $\Delta $ are independent. 

The MLE $\hat\vartheta _\varepsilon=\left( \hat\rho _\varepsilon,\hat\kappa
_\varepsilon \right) $ is consistent and it's components $\hat\rho
_\varepsilon $ and $\hat\kappa _\varepsilon $ are asymptotically independent
and have limit distributions with
different normalizing rates
\begin{align*}
&\frac{\hat\rho _\varepsilon-\rho _0}{\varepsilon
    ^{\frac{1}{H}}}\Longrightarrow  \hat \xi ,
&\frac{\hat\kappa  _\varepsilon-\kappa  _0}{\varepsilon
    }\Longrightarrow \frac{\Delta }{\II\left(\kappa _0\right)}\sim {\cal
  N}\left(0, \II\left(\kappa_0 \right)^{-1}\right).
\end{align*}
The proof follows the mains steps of the proof of the Theorem 1 is cumbersome
and do not presented here.

\bigskip

{\bf Acknowledgment.} This work was done under partial financial support of
the grant of  RSF number 14-49-10079.

\end{document}